\documentclass[12pt]{amsart}
\usepackage{amsmath,amsfonts,amsthm,amssymb,mathrsfs,
amsxtra,amscd,latexsym}
\usepackage[hmargin=2cm,vmargin=3cm]{geometry}
\usepackage[dvipdfm, usenames]{color}
 \newtheorem{thm}{Theorem}[section]

 \newtheorem{lem}[thm]{Lemma}

 \newtheorem{dfn}[thm]{Definition}

 \newtheorem{rmk}[thm]{Remark}
 \theoremstyle{definition}
 \theoremstyle{remark}
 \numberwithin{equation}{section}

\newcommand{\sm}{\left(\begin{smallmatrix}}
\newcommand{\esm}{\end{smallmatrix}\right)}
\newcommand{\mat}{\left(\begin{matrix}}
\newcommand{\emat}{\end{matrix}\right)}

\newcommand{\mb}{\mathbb}
\newcommand{\mbf}{\mathbf}

\newcommand{\mc}{\mathcal}

\def\CC{\mathbb{C}}
\def\HH{\mathbb{H}}

\def\QQ{\mathbb{Q}}

\def\ZZ{\mathbb{Z}}

\def\det{\mathrm{det}}

\def\GL{\mathrm{GL}}

\def\SL{\mathrm{SL}}

\begin{document}

\title{Cohomological relation between Jacobi forms and skew-holomorphic Jacobi forms}


 \author{Dohoon Choi}
\author{Subong Lim}

 \address{School of liberal arts and sciences, Korea Aerospace University, 200-1, Hwajeon-dong, Goyang, Gyeonggi 412-791, Republic of Korea}
  \email{choija@kau.ac.kr}	

 \address{Korea Institute for Advanced Study(KIAS), 85 Hoegiro (Cheongnyangni-dong 207-43), Dongdaemun-gu, Seoul 130-722, Republic of Korea}
  \email{subong@kias.re.kr}

 \thanks{Keynote:  Jacobi form, Eichler cohomology}

\begin{abstract}
Eichler and Zagier developed a theory of Jacobi forms to understand and extend Maass' work on the Saito-Kurokawa conjecture. Later Skoruppa introduced skew-holomorphic Jacobi forms, which play an important role in understanding liftings of modular forms and Jacobi forms. In this paper, we explain a relation between holomorphic Jacobi forms and skew-holomorphic Jacobi forms in terms of a group cohomology. More precisely, we introduce an isomorphism from the direct sum of the space of Jacobi cusp forms on $\Gamma^J$ and the space of skew-holomorphic Jacobi cusp forms on $\Gamma^J$ with the same half-integral weight to the Eichler cohomology group of $\Gamma^J$ with a coefficient module coming from polynomials.
\end{abstract}

\maketitle

\section{Introduction}

Let $\HH$ be the complex upper half plane.
A Jacobi form is a function of two variables $\tau\in\HH$ and $z\in\CC$, which satisfies modular transformation properties with respect to $\tau$ and elliptic transformation properties with respect to $z$.
Eichler and Zagier \cite{EZ} systematically developed a theory of Jacobi forms. The theory of Jacobi forms has grown enormously since then with beautiful applications in many areas of mathematics and physics, for example, the Saito-Kurokawa conjecture (see \cite{Zag}), mock theta functions (see \cite{Z}), the theory of Donaldson invariants of $\CC\mb{P}^2$ that are related to gauge theory (see \cite{GZ}) and the Mathieu moonshine (see \cite{EO}). On the other hand, skew-holomorphic Jacobi forms are introduced by Skoruppa \cite{Sko1, Sko2}. They play a crucial role in understanding liftings of modular forms and Jacobi forms.
Moreover, a Jacobi form and a skew-holomorphic Jacobi form have deep connections to modular forms through the beautiful theory of the  theta expansion, developed by Eichler and Zagier \cite{EZ},  which gives an isomorphism between  (skew-holomorphic) Jacobi forms and vector-valued modular forms.

Our aim is to explain the relation between holomorphic Jacobi forms and skew-holomorphic Jacobi forms in terms of the Eichler cohomology group of a Jacobi group.
Eichler in \cite{E} conceived the Eichler cohomology theory while studying generalized abelian integrals, which are now called the Eichler integrals, and proved that the direct sum of two spaces of cusp forms on $\Gamma$ with the same integral weight is isomorphic with the first cohomology group of $\Gamma$ with a certain module of polynomials as a coefficient if $\Gamma$ is a finitely generated Fuchsian group of the first kind which has at least one parabolic class. After the results of Eichler, the Eichler cohomology theory was studied further by Gunning \cite{G2}, Husseini and Knopp \cite{HK} and Lehner \cite{L3}.

In this paper, we deal with the Eichler cohomology theory for Jacobi forms and skew-holomorphic Jacobi forms with a half-integral weight  using a coefficient module analogous to the module of polynomials used by Eichler in \cite{E}.
More precisely, let $\Gamma\subset\SL(2,\ZZ)$ be a $H$-group. Here, an $H$-group is a finitely generated Fuchsian group of the first kind, which has at least one parabolic class. We assume that $-I\in\Gamma$ and $[\SL(2,\ZZ):\Gamma]<\infty$.
Let $k\in\ZZ+\frac12$ and $\chi$ be a (unitary) multiplier system of weight $k$ on $\Gamma$.
 We denote the vector space of Jacobi cusp forms (resp. skew-holomorphic Jacobi cusp forms) of weight $k$, index $m$ and multiplier system $\chi$ on $\Gamma^{J}$ by $S_{k,m,\chi}(\Gamma^{J})$ (resp. $J^{sk,cusp}_{k,m,\chi}(\Gamma^{J})$), where   $\Gamma^{J}$ denotes the {\it Jacobi group} $\Gamma \ltimes \ZZ^2$.

Now we consider the following special $\mathbb{C}[\Gamma^{J}]$-module to define a cohomology group of $\Gamma^{J}$.

\begin{dfn}\label{dfnofpme}
Let $P^e_{k,m}$ be the set of holomorphic functions $g(\tau,z)$ on $\HH\times\CC$ which satisfy the following conditions
\begin{enumerate}
\item[(1)] for all $X\in\ZZ^2$, $(g|_{m}X)(\tau,z) = g(\tau,z)$,
\item[(2)] $L^{k+1}_{m}g = 0$.
\end{enumerate}
\end{dfn}
This set $P_{k,m}^e$ is preserved under the slash operator $|_{-k+\frac 12,m,\chi}(\gamma,X)$ for a non-zero integer $k$ and $(\gamma,X)\in \Gamma^{J}$ and forms a vector space over $\CC$ (see section \ref{section3.1} for the detail of slash operator).
Using the space $P_{k,m}^e$ as a coefficient module, we can define a parabolic cohomology group $\tilde{H}^1_{-k+\frac12, m,\chi}(\Gamma^J,P_{k,m}^e)$ (for the precise definition see section \ref{section3.3}).
 In the following theorem we give a group cohomological relation between Jacobi cusp forms and skew-holomorphic Jacobi cusp forms of half-integral weights.
In other words,
we prove the existence of an isomorphism between the parabolic cohomology group of $\Gamma^{J}$ and the direct sum of the space of Jacobi cusp forms on $\Gamma^J$ and the space of skew-holomorphic Jacobi cusp forms on $\Gamma^J$.

\begin{thm} \label{main1} For a positive integer $k$, index $m\in \ZZ$ with $m>0$ and multiplier system $\chi$ of weight $-k+\frac 12$, we have an isomorphism
\[\tilde{\eta}: S_{k+2+\frac 12, m,\chi}(\Gamma^{J}) \oplus J^{sk,cusp}_{k+2+\frac12,m,\chi}(\Gamma^J) \cong \tilde{H}^1_{-k+\frac 12,m,\chi}(\Gamma^{J},P^e_{k,m}).\]
\end{thm}

One can be interested in finding a cohomology group isomorphic to the space of cusp forms itself.
In this direction it was studied in \cite{G2,K2}.
 Especially Knopp \cite{K} established an isomorphism between the space of cusp forms with a real weight on $\Gamma$ and the Eichler cohomology group with a special coefficient module. We would like to point out that in these cases the coefficient module is the space of  holomorphic functions satisfying a certain growth condition on which $\Gamma$ acts in a manner analogous to the action of $\Gamma$ on vector spaces of polynomials.

For Jacobi forms Choie and the second author in \cite{CL} introduced the Jacobi integrals, which are analogues of Eichler integrals. Authors of \cite{CL} also investigated their properties and gave various examples of Jacobi integrals including the construction of Jacobi integrals using the theta expansion.
 In \cite{CL2} authors studied the Eichler cohomology theory associated with Jacobi forms of arbitrary real weight using a coefficient module $\mc{P}^e_m$  of holomorphic functions on $\HH\times\CC$ satisfying a certain growth condition and proved the existence of an isomorphism between the cohomology group of $\Gamma^J$ and the space of Jacobi cusp forms on $\Gamma^J$ with a real weight.

The remainder of this paper is organized as follows. In section \ref{section2}, we introduce
vector-valued modular forms, review the supplementary function theory in terms of vector-valued Poincar\'e serie and describe how to construct the vector-valued Eichler integral for given parabolic cocycle.
In section \ref{section3},
we give basic notions of Jacobi forms and skew-holomorphic Jacobi forms and we define the cohomology group of $\Gamma^J$ with a coefficient module $P_{k,m}^e$. In section \ref{section4}, we prove Theorem \ref{main1}.

\section{Vector-valued Modular Forms} \label{section2}

In this section, we introduce the basic notions of vector-valued modular forms and vector-valued Poincar\'e series. With this vector-valued Poincar\'e series we review the supplementary function theory which was established by Knopp and his collaborators in \cite{HK, Kon}. We also describe how to construct the vector-valued Eichler integral for a given parabolic cocycle.

\subsection{Vector-valued modular forms} \label{section2.1}
We begin by introducing the definition of the vector-valued modular forms.
  Let $k\in\ZZ$ and $\chi$ be a character of $\Gamma$.
Let $p$ be a positive integer and $\rho:\Gamma\to \GL(p,\CC)$ a $p$-dimensional unitary complex representation.  We denote the standard basis elements of the vector space $\CC^p$ by $\mbf{e}_j$ for $1\leq j\leq p$. With these setups, the definition of the vector-valued modular forms are given as follows.

\begin{dfn} \label{dfnofvvm} A vector-valued weakly holomorphic modular form of weight $k$, multiplier system $\chi$ and type $\rho$ on $\Gamma$ is a sum $f(\tau) = \sum_{j=1}^p f_j(\tau)\mbf{e}_j$ of functions holomorphic in  $\HH$ satisfying the following conditions:
\begin{enumerate}
\item for all $\gamma = \sm a&b\\c&d\esm\in\Gamma$, we have
$(f|_{k,\chi,\rho} \gamma)(\tau) = f(\tau)$,

\item for each $\gamma=\sm a&b\\c&d\esm\in\SL(2,\ZZ)$,
the function $(c\tau+d)^{-k}f(\gamma\tau)$ has the Fourier expansion of the form
\[(c\tau+d)^{-k}f(\gamma\tau) = \sum_{j=1}^p\sum_{n\gg-\infty} a_{j,\gamma}(n)e^{2\pi i(n+\kappa_{j,\gamma})\tau/\lambda_{\gamma}}\mbf{e}_j,\]
where $\kappa_{j,\gamma}$ (resp. $\lambda_{\gamma}$) is a constant which depends on $j$ and $\gamma$ (resp. $\gamma$).
\end{enumerate}
\end{dfn}
Here, the slash operator $|_{k,\chi,\rho} \gamma$ is defined by
\[(f|_{k,\chi,\rho}\gamma)(\tau) = \chi(\gamma)^{-1}(c\tau+d)^{-k}\rho^{-1}(\gamma)f(\gamma\tau),\]
for $\gamma = \sm a&b\\c&d\esm \in \Gamma$, where $\gamma\tau = \frac{a\tau+b}{c\tau+d}$.
The space of all vector-valued weakly holomorphic modular forms $f(\tau)$ of weight $k$, multiplier system $\chi$ and type $\rho$ on $\Gamma$ is denoted by $M^!_{k,\chi,\rho}(\Gamma)$. There are subspaces $M_{k,\chi,\rho}(\Gamma)$ and $S_{k,\chi,\rho}(\Gamma)$ of {\it{vector-valued\ holomorphic\ modular\ forms}} and {\it{vector-valued cusp forms}}, respectively, for which we require that each $a_{j,\gamma}(n) = 0$ when $n+\kappa_{j,\gamma}$ is negative, respectively, non-positive.

\subsection{Vector-valued Eichler integrals} \label{section2.2}
 In this subsection, we introduce the vector-valued Eichler integrals.
 Let $P_k$ be the vector space of vector-valued functions $G(\tau) = \sum_{j=1}^p G_j(\tau)\mbf{e}_j$ holomorphic in $\HH$ such that $G_j(\tau)$ is a polynomial of degree at most $k$.
  It is worth mentioning that the space $P_k$ is preserved under the slash operator. With this underlying space, we introduce the vector-valued Eichler integrals.

\begin{dfn} Let $\rho$ be a $p$-dimensional representation $\rho: \Gamma\to\GL(p,\CC)$, $k$ an arbitrary nonnegative integer and $\chi$ a character of $\Gamma$.
A vector-valued Eichler integral of weight $-k$, character $\chi$ and type $\rho$ on $\Gamma$ is a vector-valued function $F(\tau)$ on $\HH$ satisfying
\[(F|_{-k,\chi,\rho}\gamma)(\tau) - F(\tau)\in P_k\]
for all $\gamma\in\Gamma$.
\end{dfn}

If we let $p_\gamma(\tau) = (F|_{-k,\chi,\rho}\gamma)(\tau) - F(\tau)$, then the vector-valued functions $p_\gamma(\tau)$ are called {\it{period functions}} of $F(\tau)$.
Then it turns out that $\{p_\gamma|\ \gamma\in\Gamma\}$ satisfies the following cocycle condition:
for $\gamma_1,\gamma_2\in\Gamma$,
\begin{equation} \label{cocyclevvm}
p_{\gamma_1\gamma_2}(\tau) = p_{\gamma_2}(\tau) + (p_{\gamma_1}|_{-k,\chi,\rho}\gamma_2)(\tau).
\end{equation}
A collection $\{p_\gamma |\ \gamma\in\Gamma\}$ of elements of $P_k$ satisfying (\ref{cocyclevvm}) is called a {\it{cocycle}} and a {\it{coboundary}} is a collection $\{p_\gamma|\ \gamma\in\Gamma\}$ such that
\[p_\gamma(\tau) = (p|_{-k,\chi,\rho}\gamma)(\tau) - p(\tau)\]
for $\gamma\in\Gamma$ with a fixed element $p(\tau) \in P_k$.
 A {\it{parabolic cocycle}} $\{p_\gamma|\ \gamma\in\Gamma\}$ is a collection of elements of $P_k$ satisfying (\ref{cocyclevvm}), in which for every parabolic class $\mc{B}$ in $\Gamma$ there exists a fixed element $Q_B(\tau)\in P_k$ such that
\[p_B(\tau) = (Q_B|_{-k,\chi,\rho}B)(\tau) - Q_B(\tau),\]
for all $B\in\mc{B}$.

A vector-valued generalized Poincar\'e series gives an example of the vector-valued Eichler integrals.
A vector-valued generalized Poincar\'e series was first defined by Lehner \cite{L3}.
 We recall the definition of a vector-valued generalized Poincar\'e series.
 Let $\{g_\gamma|\ \gamma\in\Gamma\}$ be a parabolic cocycle of elements of $P_k$ of weight $-k\in\ZZ$ with $\chi$ a multiplier system in $\Gamma$ of weight $-k$ and $\rho$ a  representation. Assume also that $g_{Q}(\tau) = 0$, where
\begin{equation} \label{lambda}
Q = \sm 1&\lambda\\ 0&1\esm\in\Gamma,\ \lambda>0,
\end{equation}
 is a generator of $\Gamma_\infty$. We define a {\it{vector-valued generalized Poincar\'e series}} as
 \begin{equation} \label{genpoin}
\Phi(\tau;r) = \sum_{j=1}^p \sum_{V= \sm a&b\\c&d\esm\in\mc{L}}\frac{(g_V)_j(\tau)}{(c\tau+d)^r}\mbf{e}_j,
\end{equation}
 where $r$ is a large positive even integer and $\mc{L}$ is any set in $\Gamma$ containing all transformations with different lower rows.

\begin{thm} \cite[Theorem 3.6]{CL2} \label{convergepoin} The generalized Poincar\'e series $\Phi(\tau;r)$, defined as in (\ref{genpoin}), converges absolutely and uniformly on compact subsets of $\HH$ for sufficiently large $r$.
\end{thm}

It was proved using a generalized Poincar\'e series $\Phi(\tau;r)$ the existence of the vector-valued Eichler integral for given period functions in \cite{CL2}.
 Following the literature on the Eichler cohomology theory, we introduce {\it{a left-finite expansion}} at each parabolic cusp for the consistency. The expansion at the cusp $i\infty$ has the form
\[F(\tau) = \sum_{j=1}^p \sum_{n\gg-\infty} a_0(n,j)e^{2\pi i(n+\kappa_{j,0})/\lambda_0}\mbf{e}_j,\]
where $Q=Q_0=\sm 1&\lambda_0\\0&1\esm,\ \lambda_0>0,$ is a generator of $\Gamma_\infty$ and
\begin{equation} \label{kappa}
\chi(Q)\rho(Q) = \sm e^{2\pi i\kappa_{1,0}}&&&&\\ &\cdot&&&\\ &&\cdot&&\\ &&&\cdot&\\ &&&&e^{2\pi i\kappa_{p,0}}\esm.
\end{equation}
Let $q_1,\cdots, q_t$ be the inequivalent parabolic cusps other than infinity. Suppose also that $Q_i = \sm *&*\\ c_i&d_i\esm,\ 1\leq i\leq t$,
is a parabolic generator, i.e., $Q_i$ is a generator of $\Gamma_j$, where $\Gamma_j$ is the cyclic subgroup of $\Gamma$ fixing $q_i,\ 1\leq i\leq t$.
To describe the expansion at a finite parabolic cusp $q_i$, choose $A_i := \sm 0&-1\\ 1&-q_i\esm$, so that $A_i$ has determinant $1$ and $A_i(q_i) = \infty$. Then the {\it{width}} of the cusp $q_i$ is a positive real number such that $A_i^{-1}\sm 1&\lambda_i\\0&1\esm A_i = Q_i$.
The expansion at the cusp $q_i$ has the form
\[F(\tau) = (\tau-q_i)^k\sum_{j=1}^p \sum_{n\gg-\infty} a_i(n,j)e^{\frac{-2\pi i(n+\kappa_{j,i})}{\lambda_i(\tau-q_i)}}\mbf{e}_j,\]
where
\begin{equation} \label{kappafinite}
\chi(Q_i)\rho(Q_i) = \sm e^{2\pi i\kappa_{1,i}}&&&&\\ &\cdot&&&\\ &&\cdot&&\\ &&&\cdot&\\ &&&&e^{2\pi i\kappa_{p,i}}\esm,
\end{equation}
for $0\leq \kappa_{j,i}<1,\ 1\leq j\leq p$ and $1\leq i\leq t$. (The motivation for the form of these expansions can be found in \cite[pp. 17-20]{K2}.)
 With these notations the following theorem gives  the existence of the vector-valued Eichler integral for given period functions.

\begin{thm} \cite[Theorem 3.9]{CL2} \label{converse}
Assume that $k$ is a nonnegative integer. Let $\chi$ be a character of $\Gamma$ and $\rho$ a unitary representation. Suppose $\{g_\gamma|\ \gamma\in\Gamma\}$ is a parabolic cocycle of weight $-k$, character $\chi$ and type $\rho$ on $\Gamma$ in $P_k$.
Then there exists a vector-valued function $\Phi(\tau)$, holomorphic in $\HH$, such that
\[(\Phi|_{-k,\chi,\rho}\gamma)(\tau) - \Phi(\tau) = g_\gamma(\tau)\]
for all $\gamma\in\Gamma$.
\end{thm}

\begin{rmk} \label{leftexpansion}
Since $\{g_\gamma|\ \gamma\in\Gamma\}$ is a parabolic cocycle, there is an element $g_i(\tau)\in P_k$ such that
\[g_{Q_i}(\tau) = (g_i|_{-k,\nu,\rho}Q_i)(\tau) - g_i(\tau),\]
for $0\leq i\leq t$.
Then $\Phi(\tau) - g_i(\tau)$ is invarinat under the slash operator $|_{-k,\chi,\rho}Q_i$ and hence has the expansions at parabolic cusps $q_i,\ 0\leq i\leq t$, of the forms
\begin{eqnarray*}
\Phi(\tau) &=& g_i(\tau) + (\tau-q_i)^k\sum_{j=1}^p \sum_{n\gg-\infty}a_i(n,j)e^{\frac{-2\pi i(n+\kappa_{j,i})}{\lambda_i(\tau-q_i)}}\mbf{e}_j,\ 1\leq i\leq t,\\
\Phi(\tau) &=& g_0(\tau) + \sum_{j=1}^p \sum_{n\gg-\infty} a_0(n,j)e^{\frac{2\pi i(n+\kappa_{j,0})}{\lambda_0}}\mbf{e}_j,\ i=0.
\end{eqnarray*}
\end{rmk}

\subsection{Vector-valued Poincar\'e series and supplementary functions} \label{section2.3}
In this subsection we review the theory of Poincar\'e series for vector-valued modular forms and introduce the supplementary function theory, following \cite{CL3}.

\begin{dfn} Fix integers $n$ and $\alpha$ with $1\leq \alpha\leq p$. The Poincar\'e series $P_{n,\alpha,\chi,\rho}(\tau)$ is defined as
\begin{equation} \label{poincare}
P_{n,\alpha,\chi,\rho}(\tau) := \frac12\sum_{\gamma=\sm a&b\\c&d\esm} \frac{e^{2\pi i(-n+\kappa_\alpha)\gamma\tau/\lambda}}{\chi(\gamma)(c\tau+d)^k}\rho(\gamma)^{-1}\mbf{e}_\alpha,
\end{equation}
where $\gamma = \sm a&b\\c&d\esm$ ranges over a set of coset representatives for $<Q>\setminus \Gamma$ and $\kappa_\alpha = \kappa_{\alpha,0}$ is as in (\ref{kappa}).
\end{dfn}

The series (\ref{poincare}) is well-defined and     invariant with respect to the action $|_{k,\chi,\rho}$ of $\Gamma$ if we assume the absolute convergence of the series $P_{n,\alpha,\chi,\rho}(\tau)$. It is known that
if $k>2$ then the component function $(P_{n,\alpha,\chi,\rho})_j(\tau)$ of $P_{n,\alpha,\chi,\rho}(\tau),\ 1\leq j\leq p$, converges absolutely uniformly on compact subsets of $\HH$. In paricular, each $(P_{n,\alpha,\chi,\rho})_j(\tau)$ is holomorphic in $\HH$ (see Proposition 2.5 in \cite{CL3}) and its Fourier expansion at $i\infty$ is of the form (see \cite[Theorem 2.6]{CL3})
\begin{equation} \label{Fourierofvvpoincare}
(P_{n,\alpha,\chi,\rho})_j(\tau) = \sum_{j=1}^p \delta_{\alpha,j} e^{2\pi i(-n+\kappa_\alpha)\tau/\lambda} + \sum_{j=1}^p \sum_{l+\kappa_j>0} a(l,j)e^{2\pi i(l+\kappa_j)\tau/\lambda}.
\end{equation}

Now we look at the properties of the Poincar\'e series $P_{n,\alpha,\chi,\rho}(\tau)$.
One can see that the function $P_{n,\alpha,\chi,\rho}(\tau)$ vanishes at all cusps of $\Gamma$ which are not equivalent to $i\infty$ and that $S_{k,\chi,\rho}(\Gamma)$ is spanned by Poincar\'e series $P_{n,\alpha,\chi,\rho}(\tau)$ with $n+\kappa_\alpha>0$ (see Theorem 2.7 in \cite{CL3}).

Suppose $f(\tau)\in S_{k+2,\chi,\rho}(\Gamma)$ with $k>0$. Then there exist complex numbers $b_1,\cdots, b_s$ such that $f(\tau) = \sum_{i=1}^s b_i P_{n_i,\alpha_i,\chi,\rho}(\tau)$. Put $f^*(\tau) = \sum_{i=1}^s\overline{b_i}P_{n'_i,\alpha_i,\bar{\chi},\bar{\rho}}(\tau)$, where
\begin{equation*}
n'_i =
\begin{cases}
-n_i & \text{if $\kappa_\alpha =0$},\\
1 - n_i & \text{if $\kappa_\alpha>0$}.
\end{cases}
\end{equation*}
Note that
if we let
\[\bar{\chi}(Q)\bar{\rho}(Q) = \sm e^{2\pi i\kappa'_1}& & & & \\
						                 & \cdot &&&\\
						                 && \cdot &&\\
						                 &&& \cdot &\\
						                 &&&& e^{2\pi i\kappa'_p}\esm,\]
with $0\leq\kappa_j'<1$ for $1\leq j\leq p$, then we have
\begin{equation*}
\kappa_j' =
\begin{cases}
0 & \text{if $\kappa_j=0$},\\
1-\kappa & \text{if $\kappa_j>0$},
\end{cases}
\end{equation*}
for $1\leq j\leq p$.
Thus we have the expansion at $i\infty$
\begin{eqnarray*}
P_{n'_i,\alpha_i,\bar{\chi},\bar{\rho}}(\tau) &=& e^{2\pi i(-n'_i+\kappa'_{\alpha_i})\tau/\lambda}\mbf{e}_{\alpha_i}+\sum_{j=1}^p \sum_{l+\kappa_j>0} a_{n'_i,\alpha_i}(l,j)e^{2\pi i(l+\kappa_j)\tau/\lambda}\mbf{e}_j\\
&=& e^{2\pi i(n_i-\kappa_{\alpha_i})\tau/\lambda}\mbf{e}_{\alpha_i} + \sum_{j=1}^p \sum_{l+\kappa_j>0} a_{n'_i,\alpha_i}(l,j)e^{2\pi i(l+\kappa_j)\tau/\lambda}\mbf{e}_j.
\end{eqnarray*}
It follows that $f^*(\tau)\in M^!_{k+2,\bar{\chi},\bar{\rho}}(\Gamma)$, $f^*(\tau)$ has a pole at $i\infty$ with principal part
\[\sum_{i=1}^s \overline{b_i}e^{2\pi i(n_i-\kappa_{\alpha_i})\tau/\lambda},\]
and $f^*(\tau)$ vanishes at all of the other cusps of $\Gamma$. We call $f^*(\tau)$ the {\it{function supplementary to}} $f(\tau)$.

Functions $f(\tau)$ and $f^*(\tau)$ have important relations which can be expressed in terms of period functions.  A form $f(\tau)\in M^!_{k+2,\chi,\rho}(\Gamma)$ is a {\it{vector-valued weakly holomorphic cusp form}} if its constant term vanishes. Let $S_{k+2,\chi,\rho}^!(\Gamma)$ denote the space of vector-valued weakly holomorphic cusp forms.

Suppose that $f(\tau)\in S_{k+2,\chi,\rho}(\Gamma)$. Then we have a nonholomorphic  Eichler integral of $f(\tau)$
\begin{eqnarray*}
\mc{E}^N_f(\tau) &:=& \frac{1}{c_{k+2}}\biggl[\int^{i\infty}_\tau f(z)(\bar{\tau}-z)^{k}dz\biggr]^-,
\end{eqnarray*}
where $c_k = -\frac{(k-2)!}{(2\pi i)^{k-1}}$ and $[\ ]^-$ indicates the complex conjugate of the function inside $[\ ]^-$.
On the other hand, suppose that $f(\tau) =\displaystyle  \sum_{j=1}^p\sum_{n\gg-\infty\atop n+\kappa_j\neq0} a(n,j)e^{2\pi i(n+\kappa_j)\tau/\lambda}\mbf{e}_j \in S^!_{k+2,\chi,\rho}(\Gamma)$. Then we have a holomorphic Eichler integral of $f(\tau)$
\begin{eqnarray*}
\mc{E}^H_f(\tau) &:=&  \sum_{j=1}^p\sum_{n\gg-\infty\atop n+\kappa_j\neq0} a(n,j)\biggr(\frac{n+\kappa_j}{\lambda}\biggl)^{-(k+1)}e^{2\pi i(n+\kappa_j)\tau/\lambda}\mbf{e}_j + c_f.
\end{eqnarray*}
Here, the constant term $c_f$ is determined by the fact that the Fourier coefficients of modular forms of negative weight are completely determined by the principal part of the expansion of those forms at the cusps (see \cite{Leh}).
For example, if we assume that $f(\tau)$ has a pole at $i\infty$ and that it is holomorphic at all other cusps, then $c_f$ is equal to
\[c_f := \sum_{j=1}^p \delta_{\kappa_j,0}\biggl(\frac{1}{\lambda(k+1)!}\sum_{t=1}^p\sum_{l<0}\sum_{\gamma=\sm a&b\\c&d\esm\in C^+}a(l,t) \biggl(\frac{-2\pi i}{c}\biggr)^{k+2}\chi^{-1}(\gamma)\rho(\gamma^{-1})_{j,t}e^{\frac{2\pi i}{c\lambda}(l+\kappa_t)a}\biggr)\mbf{e}_j.\]

We introduce the {\it{period functions}} for $f(\tau)$ by
\begin{eqnarray*}
r^H(f,\gamma;\tau) &:=& c_{k+2}(\mc{E}^H_f - \mc{E}^H_f|_{-k,\chi,\rho}\gamma)(\tau),\\
r^N(f,\gamma;\tau) &:=& c_{k+2}(\mc{E}^N_f - \mc{E}^N_f|_{-k,\bar{\chi},\bar{\rho}}\gamma)(\tau),
\end{eqnarray*}
where $\gamma\in\Gamma$.
For the period functions we have the following theorem.

\begin{thm} \cite[Theorem 2.8]{CL3}, \cite[Section 2]{HK} \label{suppleperiod} Suppose that $k$ is a positive integer and $f(\tau)\in S_{k+2,\chi,\rho}(\Gamma)$.
Then we have
\begin{enumerate}
\item  $r^H(f,\gamma;\tau) = [r^H(f^*,\gamma;\bar{\tau})]^-$ for all $\gamma\in\Gamma$,
\item  $r^H(f,\gamma;\tau) = [r^N(f,\gamma;\bar{\tau})]^-$ for all $\gamma\in\Gamma$.
\end{enumerate}
\end{thm}

\section{Jacobi forms and associated cohomology groups } \label{section3}

In this section, we review basic notions of Jacobi forms (see \cite{EZ, Z0}) and  skew-holomorphic Jacobi forms (see \cite{Sko1,Sko2}). We also review the definition of Jacobi integrals (see \cite{CL}) and define the cohomology group for Jacobi forms of half-integral weight.

\subsection{Jacobi forms} \label{section3.1}
First we fix some notations.
Let $\Gamma^{J}=\Gamma\ltimes \ZZ^2$ be a {\it{Jacobi group}} with associated composition law
\[(\gamma_1,(\lambda_1,\mu_1))\cdot(\gamma_2,(\lambda_2,\mu_2))=(\gamma_1 \gamma_2, (\widetilde{\lambda_1}+\lambda_2,\widetilde{\mu_1}+\mu_2)),\]
where $(\widetilde{\lambda},\widetilde{\mu})=(\lambda,\mu)\cdot \gamma_2$.
Then $\Gamma^{J}$ acts on $\HH\times\CC$ as a group of automorphism. The action is given by
\[(\gamma, (\lambda,\mu))\cdot(\tau,z)= \biggl(\gamma\tau,\frac{z+\lambda\tau+\mu}{c\tau+d}\biggr),\]
where $\gamma\tau = \frac{a\tau+b}{c\tau+d}$ for $\gamma = \sm a&b\\c&d\esm\in\Gamma$.
Let $k$ be a half-integer and $\chi$ be a multiplier system of weight $k$ on $\Gamma$, i.e., $\chi : \Gamma\to\CC$ satisfies
\begin{enumerate}
\item[(1)] $|\chi(\gamma)| = 1$ for all $\gamma\in\Gamma$,
\item[(2)] $\chi$ satisfies the consistency condition
\[\chi(\gamma_3)(c_3\tau+d_3)^{k}=\chi(\gamma_1)\chi(\gamma_2)(c_1\gamma_2\tau+d_1)^{k}(c_2\tau+d_2)^{k},\]
where $\gamma_3=\gamma_1\gamma_2$ and $\gamma_i =\sm *&*\\c_i&d_i\esm, i =1,2$ and $3$,
\item[(3)] $\chi$ satisfies  $\chi(-I) = e^{\pi ik}$.
\end{enumerate}
For $\gamma=\sm a&b\\c&d\esm \in\Gamma , X = (\lambda, \mu)\in \ZZ^2$ and $m\in\ZZ$ with $m>0$, we define
\[(\Phi|_{k,m,\chi} \gamma)(\tau,z) := (c\tau+d)^{-k}\bar{\chi}(\gamma)e^{-2\pi im\frac{cz^2}{c\tau+d}}\Phi(\gamma(\tau,z))\]
and
\[(\Phi|_{m}X)(\tau,z) :=e^{2\pi im(\lambda^2\tau+ 2\lambda z+\mu\lambda)}\Phi(\tau,z+\lambda\tau+\mu),\]
where $\gamma(\tau,z) = (\frac{a\tau+b}{c\tau+d},\frac{z}{c\tau+d})$.
Then $\Gamma^{J}$ acts on the space of  functions on $\HH\times\CC$ by
\begin{equation} \label{slashoperator}
(\Phi|_{k,m,\chi} (\gamma,X))(\tau,z) := (\Phi|_{k,m,\chi} \gamma |_{m}X)(\tau,z).
\end{equation}

We introduce the definition of a Jacobi form.

\begin{dfn}
A {\it{Jacobi form}} of weight $k$, index $m$ and multiplier system $\chi$ on $\Gamma^J$  is a holomorphic mapping $\Phi(\tau,z)$ on $\HH\times\CC$ satisfying
\begin{enumerate}
\item[(1)] $(\Phi|_{k,m,\chi} \gamma)(\tau,z) =\Phi(\tau,z)$ for every $\gamma\in\Gamma$,
\item[(2)] $(\Phi|_{m}X)(\tau,z) = \Phi(\tau,z)$ for every $X\in \ZZ^2$,
\item[(3)] for each $\gamma=\sm a&b\\c&d\esm \in\SL(2,\ZZ)$, the function $(c\tau+d)^{-k}e^{2\pi im\frac{-cz^2}{c\tau+d}}\Phi((\gamma,0)\cdot(\tau,z))$ has the Fourier expansion of the form
\begin{equation} \label{Jacobifourier}
\sum_{l,r\in\ZZ\atop 4(l+\kappa_\gamma)-\lambda_\gamma mr^2 \geq 0}a(l,r)e^{2\pi i(l+\kappa_{\gamma})/\lambda_\gamma}e^{2\pi irz},
\end{equation}
with a suitable $0\leq \kappa_\gamma<1,\ \lambda_\gamma\in\ZZ$.
\end{enumerate}
If a  Jacobi form satisfies the  condition $a(l,r)\neq 0$ only if 
$4(l+\kappa_\gamma)-\lambda_\gamma mr^2 >0$, then it is called a {\it {Jacobi cusp form}}.
\end{dfn}

Now we look into the theta series, which plays an important role in the proofs of our main theorems.
Let $S$ be a positive integer  and $a,b\in\QQ$. We consider the theta series
\[\theta_{S,a,b}(\tau,z) := \sum_{\lambda\in\ZZ}e^{\pi iS((\lambda+a)^2\tau+2(\lambda+a)(z+b))}\]
with characteristic $(a,b)$ converging normally on $\HH\times\CC$.

We are in a position to explain  the beautiful theory of the  theta expansion, developed by Eichler and Zagier \cite{EZ},  which gives an isomorphism between  (skew-holomorphic) Jacobi forms and vector-valued modular forms.

\begin{thm} \cite[Section 3]{Z0} \label{decomposition} Let $\Phi(\tau,z)$ be holomorphic as a function of $z$ and satisfy
\begin{equation} \label{elliptic}
(\Phi|_{m} X)(\tau,z) = \Phi(\tau,z)\ \text{for every}\ X\in \ZZ^2.
\end{equation}
Then we have
\begin{equation} \label{thetaexpansion}
\Phi(\tau,z) = \sum_{a\in\mathcal{N}}f_a(\tau)\theta_{2m,a,0}(\tau,z)
\end{equation}
with uniquely determined holomorphic functions $f_a:\HH\to\CC$, where $\mc{N} = \ZZ/2m\ZZ$.
If $\Phi(\tau,z)$ also satisfies the transformation
\[(\Phi|_{k,m,\chi} \gamma)(\tau,z) =\Phi(\tau,z)\ \text{for every}\ \gamma\in\SL(2,\ZZ),\]
then we have for each $a\in\mc{N}$
\[f_a\biggl(-\frac1\tau\biggr) = \chi(\sm 0&-1\\1&0\esm)\det\biggl(\frac \tau i\biggr)^{-\frac j2}\tau^k(2m)^{-\frac 12}\sum_{b\in\mathcal{N}}e^{2\pi i(2mab)}f_b(\tau)\]
and
\[f_a(\tau+1) = \chi(\sm 1&1\\0&1\esm)e^{-2\pi ima^2}f_a(\tau).\]
Furthermore, if $\Phi(\tau,z)$ is a Jacobi form in $J_{k,m,\chi}(\Gamma^{(1,j)})$, then
functions in $\{f_a|\ a\in\mathcal{N}\}$ necessarily must have the Fourier expansions of the form
\[f_a(\tau) =\sum_{l\geq0\atop \text{rational}} a(l)e^{2\pi il\tau}.\]
\end{thm}

The decomposition by theta functions as in (\ref{thetaexpansion}) is called the {\it{theta expansion}}.
Now we explain the isomorphism between Jacobi forms and vector-valued modular forms induced by the theta expansion more precisely. First, we can define the space of vector-valued modular forms associated with $J_{k,m,\chi}(\Gamma^{J})$ by the theta expansion.
Let $\Phi(\tau,z)$ be a Jacobi form in $J_{k,m,\chi}(\Gamma^{J})$. By Theorem \ref{decomposition}, we have the theta expansion
\[\Phi(\tau,z) = \sum_{a\in\mc{N}}f_a(\tau)\theta_{2m,a,0}(\tau,z).\]
We take a multiplier system $\chi''$ of weight $1/2$ on $\SL(2,\ZZ)$, for example we can take a power of eta-multiplier system: $\chi''(\gamma) = \frac{\eta(\gamma\tau)}{\eta(\tau)}$ for $\gamma\in\SL(2,\ZZ)$, where $\eta(\tau) = e^{\frac{\pi i\tau}{12}}\prod_{n=1}^\infty (1-e^{2\pi in\tau})$ is the Dedekind eta function.
Then we define a representation $\rho' : \SL(2,\ZZ) \to \GL(|\mc{N}|,\CC)$ by
\begin{equation} \label{doublerho1}
\rho'(T)\mbf{e}_a = \chi''(T)e^{-2\pi ima^2}\mbf{e}_a
\end{equation}
and
\begin{equation} \label{doublerho2}
\rho'(S)\mbf{e}_a = \chi''(S)\frac{i^{\frac 12}}{\sqrt{2m}}\sum_{b\in\mc{N}}e^{2\pi i(2mab)}\mbf{e}_b,
\end{equation}
where $T = \sm 1&1\\0&1\esm$ and $S = \sm 0&-1\\1&0\esm$.
Let $\chi'$ be a character of $\Gamma$ defined by
\begin{equation} \label{doublechi}
\chi'(\gamma) = \chi(\gamma)\overline{\chi''}(\gamma).
\end{equation}
Then a vector-valued function $\sum_{a\in\mc{N}} f_a(\tau)\mbf{e}_a$ is a vector-valued modular form in  $M_{k-\frac 12, \chi', \rho'}(\Gamma)$.

\begin{thm} \cite[Theorem 3.3]{Z0} \label{isomorphism}
The theta expansion  gives an isomorphism between $J_{k,m,\chi}(\Gamma^{J})$ and $M_{k-\frac 12, \chi', \rho'}(\Gamma)$. Furthermore, this isomorphism sends Jacobi cusp forms to vector-valued cusp forms.
\end{thm}

\subsection{Skew-holomorphic Jacobi forms} \label{section3.2}
Skoruppa \cite{Sko1, Sko2} introduced skew-holomorphic Jacobi forms, which play a crucial role in understanding liftings of modular forms and Jacobi forms. Now we briefly discuss the definition of skew-holomorphic Jacobi forms. For a half-integer $k$ and a integer $m$, we have the following slash operator on functions $\phi:\HH\times\CC\to\CC$:
\[(\phi|_{k,m,\chi}^{sk} (\gamma,X))(\tau,z) := \bar{\chi}(\gamma)(c\bar{\tau}+d)^{1-k}|c\tau+d|^{-1}e^{2\pi im(-\frac{c(z+\lambda\tau+\mu)^2}{c\tau+d}+\lambda^2\tau+2\lambda z)}\phi((\gamma,X)\cdot(\tau,z)),\]
for all $(\gamma,X) = \biggl(\sm a&b\\c&d\esm, (\lambda,\mu)\biggr)\in\Gamma^{J}$, where $\chi$ is a multiplier system of weight $k$ on $\Gamma$.

\begin{dfn} A function $\phi:\HH\times\CC\to\CC$ is a skew-holomorphic Jacobi form of weight $k$, index $m$ and multiplier system $\chi$ on $\Gamma^{J}$ if $\phi(\tau,z)$ is real-analytic in $\tau\in\HH$, is holomorphic in $z\in\CC$, and satisfies the following conditions:
\begin{enumerate}
\item[(1)] for all $(\gamma,X)\in\Gamma^{J}$, $(\phi|_{k,m,\chi}^{sk}(\gamma,X))(\tau,z) = \phi(\tau,z)$,
\item[(2)] for each $\gamma=\sm a&b\\c&d\esm \in\SL(2,\ZZ)$, the function $(c\tau+d)^{-k}e^{2\pi im(-\frac{cz^2}{c\tau+d})}\Phi((\gamma,0)\cdot(\tau,z))$ has the Fourier expansion of the form
\begin{equation} \label{Jacobifourier}
\sum_{l,r\in\ZZ\atop 4m(l+\kappa_\gamma)-\lambda_\gamma r^2\leq 0}a(l,r)e^{\frac{-\pi(\lambda_\gamma r^2-4m(l+\kappa_{\gamma}))v}{m\lambda_\gamma}}e^{2\pi i\tau(l+\kappa_{\gamma})/\lambda_\gamma}e^{2\pi irz},
\end{equation}
with a suitable $0\leq \kappa_\gamma<1,\ \lambda_\gamma\in\ZZ$.
\end{enumerate}
\end{dfn}

If the Fourier expansion in (\ref{Jacobifourier}) is only over $4m(l+\kappa_\gamma)-\lambda_\gamma r^2<0$,  then $\phi(\tau,z)$ is a {\it{skew-holomorphic Jacobi cusp form}} of weight $k$, index $m$ and multiplier system $\chi$ on $\Gamma^{J}$. We denote the spaces of  skew-holomorphic Jacobi forms and skew-holomorphic Jacobi cusp forms, each of weight $k$, index $m$ and multiplier system $\chi$ on $\Gamma^{J}$, by $J^{sk}_{k,m,\chi}(\Gamma^{J})$ and $J^{sk,cusp}_{k,m,\chi}(\Gamma^{J})$, respectively.

Bringmann and Richter in \cite{BR} found the close relationship between skew-holomorphic Jacobi forms and vector-valued modular forms.
Like a Jacobi form, a skew-holomorphic Jacobi form has the theta decomposition. More precisely, if $\phi(\tau,z)\in J^{sk}_{k,m,\chi}(\Gamma^{J})$, then we have
\[\phi(\tau,z) = \sum_{a\in \mc{N}}f_a(\tau)\theta_{2m,a,0}(\tau,z).\]
But in this case, $\sum_{a\in\mc{N}}\overline{f_a(\tau)}\mbf{e}_a$ is a vector-valued modular form in $M_{k-\frac12, \overline{\chi'},\overline{\rho'}}(\Gamma)$ where $\chi'$ and $\rho'$ are defined in (\ref{doublechi}) and (\ref{doublerho1}), (\ref{doublerho2}), respectively.
This gives an isomorphism between skew-holomorphic Jacobi forms and vector-valued modular forms.

\begin{thm} \cite[Section 6]{BR}  \label{isomorphismskew}
The theta expansion  gives an isomorphism between $J^{sk}_{k,m,\chi}(\Gamma^{J})$ and $M_{k-\frac 12, \overline{\chi'}, \overline{\rho'}}(\Gamma)$. Furthermore, this isomorphism sends skew-holomorphic Jacobi cusp forms to vector-valued cusp forms.
\end{thm}

\subsection{Cohomology groups} \label{section3.3}
We give the definition of the coefficient module in Definition \ref{dfnofpme} in terms of the kernel of the heat operator. Using an analogy of Bol's identity for Jacobi forms and the theta expansion, we give the more explicit description of the coefficient module.

\begin{lem} \label{integraldecomposition}
Let $g(\tau,z)$ be a holomorphic function on $\HH\times\CC$. Then a function $g(\tau,z)$ is an element of $P_{k,m}^e$ if and only if a function $g(\tau,z)$ can be written as
\[g(\tau,z) = \sum_{a\in\mc{N}} g_a(\tau)\theta_{2m,a,0}(\tau,z),\]
where $g_a(\tau)$ is a polynomial of degree at most $k$ for $a\in\mc{N}$.
\end{lem}

To prove Lemma \ref{integraldecomposition} we need the following theorem, which is an analogy of Bol's identity (see \cite{B}) for the functions on $\HH\times\CC$.

\begin{thm} \cite[Theorem 3.4]{CK}
Let $\phi:\HH\times\CC\to \CC$ be a smooth function. For a nonnegative integer $k$ and $(\gamma,X)\in \Gamma^{J}$ we have
\[ \{(L_{m}^{k+1}\phi)|_{k+\frac 12+2, m,\chi}(\gamma,X)\}(\tau,z) = L_{m}^{k+1}(\phi|_{-k+\frac 12,m,\chi}(\gamma,X))(\tau,z).\]
\end{thm}

Now we give the proof of Lemma \ref{integraldecomposition}.

\begin{proof} [\bf Proof of Lemma \ref{integraldecomposition}]
Suppose that $g(\tau,z)$ is an element of $P_{k,m}^e$. Then, by Theorem \ref{decomposition}, a function $g(\tau,z)$ can be written as
\[g(\tau,z) = \sum_{a\in\mc{N}} g_a(\tau)\theta_{2m,a,0}(\tau,z),\]
where $g_a(\tau)$ is a holomorphic function on $\HH$. Note that the theta function $\theta_{2m,a,0}(\tau,z)$ is a heat kernel, i.e, $L_{m}(\theta_{2m,a,0}) = 0$. Therefore, we see that
\[L^{k+1}_{m}(g)(\tau,z) = (8\pi im)^{k+1}\sum_{a\in\mc{N}}\biggl(\biggl(\frac{\partial}{\partial\tau}\biggr)^{k+1}g_a\biggr)(\tau)\theta_{2m,a,0}(\tau,z)=0.\]
Therefore, $(\frac{\partial}{\partial\tau})^{k+1}g_a$ should be zero for every $a\in\mc{N}$. This implies that $g_a(\tau)$ is a polynomial of degree at most $k$ for $a\in\mc{N}$.

Conversely, suppose that a function $g(\tau,z)$ is given by
\[g(\tau,z) = \sum_{a\in\mc{N}} g_a(\tau)\theta_{2m,a,0}(\tau,z),\]
where $g_a(\tau)$ is a polynomial of degree at most $k$. Then one can check that a function $g(\tau,z)$ satisfies two conditions in the definition of $P_{k,m}^e$ by the similar argument as above.
\end{proof}

Note that if $\Phi(\tau,z)$ is a Jacobi integral of weight $-k+\frac 12$, index $m$ and multiplier system $\chi$ on $\Gamma^{J}$, then it satisfies
\[(\Phi|_{-k+\frac 12,m, \chi}(\gamma,X))(\tau,z) - \Phi(\tau,z) \in P_{k,m}^e\]
for every $(\gamma,X)\in\Gamma^{J}$. If we let $p_{(\gamma,X)}(\tau,z) := (\Phi|_{-k+\frac 12,m, \chi}(\gamma,X))(\tau,z) - \Phi(\tau,z)$, then a function $p_{(\gamma,X)}(\tau,z)$ is called a {\it{period function}} of $\Phi(\tau,z)$ for $(\gamma,X)\in\Gamma^{J}$. Note that a collection $\{p_{(\gamma,X)}|\ (\gamma,X)\in \Gamma^{J}\}$ satisfies a cocycle condition
\begin{equation} \label{cocycle}
p_{(\gamma_1,X_1)(\gamma_2,X_2)}(\tau,z) = (p_{(\gamma_1,X_1)}|_{-k+\frac 12,m,\chi}(\gamma_2,X_2))(\tau,z) + p_{(\gamma_2,X_2)}(\tau,z)
\end{equation}
for $(\gamma_1,X_1),\ (\gamma_2,X_2)\in \Gamma^{J}$.

The cocycle relation (\ref{cocycle}) gives a definition of the cohomology group of $\Gamma^{J}$ as follows.
A collection $\{p_{(\gamma,X)}|\ (\gamma,X)\in \Gamma^{J}\}$ of elements of $P^e_{k,m}$ satisfying
(\ref{cocycle})
is called a {\it{cocycle}} and a {\it{coboundary}} is a collection $\{p_{(\gamma,X)}|\ (\gamma,X)\in \Gamma^{J}\}$ such that
\[p_{(\gamma,X)}(\tau,z) = (p|_{-k+\frac12, m,\chi}(\gamma,X))(\tau,z) - p(\tau,z)\]
for all $(\gamma,X)\in \Gamma^{J}$ with a fixed element $p(\tau,z)$ of $P^e_{k,m}$. The {\it{cohomology group}}
$H^1_{-k+\frac 12,m,\chi}(\Gamma^{J},P^e_{k,m})$
is the quotient of the cocycles by the coboundaries. A cocycle $\{p_{(\gamma,X)}|\ (\gamma,X)\in \Gamma^{J}\}$ is called the {\it{parabolic cocycle}} if for every parabolic element $B$ in $\Gamma$ there exists a function $Q_{B}(\tau,z)$ in $P^e_{k,m}$ such that
\[p_{(B,0)}(\tau,z) = (Q_{B}|_{-k+\frac 12,m,\chi})(B,0)(\tau,z) - Q_{B}(\tau,z).\]
 The {\it{parabolic cohomology group}} $\tilde{H}^1_{-k+\frac 12,m,\chi}(\Gamma^{J},P^e_{k,m})$ is defined as the vector space obtained by forming the quotient of the parabolic cocycles by the coboundaries.

\section{Proof of Theorem \ref{main1}} \label{section4}

First, we construct the mapping $\tilde{\eta}$ from $S_{k+2+\frac 12, m,\chi}(\Gamma^{J}) \oplus J^{sk,cusp}_{k+2+\frac12,m,\chi}(\Gamma^J)$ to $\tilde{H}^1_{-k+\frac 12,m,\chi}(\Gamma^{J},P^e_{k,m})$. If we have a Jacobi cusp form $\Phi(\tau,z)\in S_{k+2+\frac12,m,\chi}(\Gamma^J)$, then we have the theta expansion
\begin{equation} \label{Jacobithetaexpansion}
\Phi(\tau,z) = \sum_{a\in\mc{N}} f_a(\tau)\theta_{2m,a,0}(\tau,z).
\end{equation}
Then $f(\tau):=\sum_{a\in\mc{N}}f_a(\tau)\mbf{e}_a$ is a vector-valued cusp form in $S_{k+2,\chi',\rho'}(\Gamma)$. The collection
\[\{r^H(f,\gamma;\tau)|\ \gamma\in\Gamma\}\]
is a parabolic cocycle in $P_{k}$, where $r^H(f,\gamma;\tau)$ is a period function of $f(\tau)$ defined in section \ref{section2.3}. Then this cocycle is parabolic (see \cite[pp. 158-159]{L3} or \cite[pp. 129]{L2}). Using this cocycle, we define a cocycle $\{r(\Phi,(\gamma,X);\tau,z)|\ (\gamma,X)\in\Gamma^J\}$ in $P^e_{k,m}$ by
\[r(\Phi,(\gamma,X);\tau,z) = \sum_{a\in\mc{N}} r^H(f,\gamma;\tau)_a\theta_{2m,a,0}(\tau,z),\]
where $r^H(f,\gamma;\tau)_a$ is the $a$th component function of a vector-valued function $r^H(f,\gamma;\tau)$.  Then one can check that this cocycle is also parabolic by using the property of the theta expansion: for $(\gamma,X)\in\Gamma^{J}$, we have
\begin{equation} \label{propertyofthetaexpansion}
\biggl(\biggl(\sum_{a\in\mc{N}}f_a\theta_{2m,a,0}\biggr)\biggr|_{-k+\frac 12, m,\chi}(\gamma,X)\biggr)(\tau,z) = \sum_{a\in\mc{N}}F_a(\tau)\theta_{2m,a,0}(\tau,z),
\end{equation}
where $\sum_{a\in\mc{N}}F_a(\tau)\mbf{e}_a = \biggl(\biggl(\sum_{a\in\mc{N}}f_a(\tau)\mbf{e}_a\biggr)\biggr|_{-k,\chi',\rho'}\gamma\biggr)(\tau)$.  We define a map $\tilde{\beta}$ from $S_{k+2+\frac 12, m,\chi}(\Gamma^{J})$ to  $\tilde{H}^1_{-k+\frac 12,m,\chi}(\Gamma^{J},P^e_{k,m})$ by
\[\tilde{\beta}(\Phi) = < r(\Phi,(\gamma,X);\tau,z)|\ (\gamma,X)\in\Gamma^J>,\]
where $< r(\Phi,(\gamma,X);\tau,z)|\ (\gamma,X)\in\Gamma^J>$ is a cocycle class represented by a parabolic cocycle $\{r(\Phi,(\gamma,X);\tau,z)|\ (\gamma,X)\in\Gamma^J\}$. Likewise, if we have a skew-Jacobi cusp form $\Psi(\tau,z)\in J^{sk,cusp}_{k+2+\frac12,m,\chi}(\Gamma^J)$, then we have the theta expansion
\begin{equation} \label{skewJacobithetaexpansion}
\Psi(\tau,z) = \sum_{a\in\mc{N}} g_a(\tau)\theta_{2m,a,0}(\tau,z).
\end{equation}
Then $g(\tau):=\sum_{a\in\mc{N}}\overline{g_a(\tau)}\mbf{e}_a$ is a vector-valued cusp form in $S_{k+2,\overline{\chi'},\overline{\rho'}}(\Gamma)$.
We define a map $\tilde{\alpha}$ from $J^{sk,cusp}_{k+2+\frac12,m,\chi}(\Gamma^J)$ to  $\tilde{H}^1_{-k+\frac 12,m,\chi}(\Gamma^{J},P^e_{k,m})$ by
\[\tilde{\alpha}(\Psi) = < r^*(\Psi,(\gamma,X);\tau,z)|\ (\gamma,X)\in\Gamma^J>,\]
where $r^*(\Phi,(\gamma,X);\tau,z) = \sum_{a\in\mc{N}} \overline{r^H(g^*, \gamma;\tau)_a}\theta_{2m,a,0}(\tau,z)$, where $g^*(\tau)$ is a function supplementary to $g(\tau)$. Then this cocycle is also parabolic and a map $\tilde{\alpha}$ is well-defined. Combining $\tilde{\alpha}$ and $\tilde{\beta}$, we define a map $\tilde{\eta} : S_{k+2+\frac 12, m,\chi}(\Gamma^{J}) \oplus J^{sk,cusp}_{k+2+\frac12,m,\chi}(\Gamma^J) \to \tilde{H}^1_{-k+\frac 12,m,\chi}(\Gamma^{J},P^e_{k,m})$ by
\[\tilde{\eta}(\Phi,\Psi) = \tilde{\beta}(\Phi) + \tilde{\alpha}(\Psi).\]

Next, we prove that a map $\tilde{\eta}$ is injective. This proof is based on the argument in \cite{KLR}.
Suppose that $(\Phi,\Psi) \in S_{k+2+\frac 12, m,\chi}(\Gamma^{J}) \oplus J^{sk,cusp}_{k+2+\frac12,m,\chi}(\Gamma^J)$ and that a cocycle $\{ r(\Phi,(\gamma,X);\tau,z) + r^*(\Phi,(\gamma,X);\tau,z)|\ (\gamma,X)\in\Gamma^J\}$ is a coboundary in $P^e_{k,m}$. Then by (\ref{propertyofthetaexpansion}) a cocycle $\{r^H(f,\gamma;\tau) + r^H(g^*,\gamma;\tau)|\ \gamma\in\Gamma\}$ is a coboundary in $P_k$. By Theorem \ref{suppleperiod}  we see that $r^H(g^*,\gamma;\tau) = r^N(g,\gamma;\tau)$. Now we use the integral representations of $r^H(f,\gamma;\tau)$ and $r^N(g,\gamma;\tau)$ (see \cite{KM2})
\begin{eqnarray*}
r^H(f,\gamma;\tau) &=& \int^{i\infty}_{\gamma^{-1}(i\infty)}f(z)(\tau-z)^kdz,\\
r^N(g,\gamma;\tau) &=& \biggl[\int^{i\infty}_{\gamma^{-1}(i\infty)}g(z)(\bar{\tau}-z)^kdz\biggr]^-.
\end{eqnarray*}
Since $\{r^H(f,\gamma;\tau) + r^N(g,\gamma;\tau)|\ \gamma\in\Gamma\}$ is a coboundary in $P_k$, there exists $p(\tau)$ in $P_k$ such that
\[r^H(f,\gamma;\tau) + r^N(g,\gamma;\tau) = (p|_{-k,\chi',\rho'}\gamma)(\tau) - p(\tau),\ \gamma\in\Gamma.\]
This implies that $(F+\hat{G}-p)(\tau)$ is invariant under the slash operator $|_{-k,\chi',\rho'}$, where $F(\tau)$ and $\hat{G}(\tau)$ are the vector-valued Eichler integrals defined by
\begin{eqnarray*}
F(\tau) &=& \int^{\tau}_{i\infty}f(z)(\tau-z)^kdz,\\
\hat{G}(\tau) &=& \biggl[\int^{\tau}_{i\infty}g(z)(\bar{\tau}-z)^kdz\biggr]^-.
\end{eqnarray*}
Now, since $g(\tau)$ is holomorphic in $\HH$, we have
\[\frac{\partial}{\partial \bar{\tau}} (g_a(\tau)\hat{G}_a(\tau)) = g_a(\tau)\frac{\partial \hat{G}_a(\tau)}{\partial\bar{\tau}} = |g_a(\tau)|^2(\tau-\bar{\tau})^k\]
for each $a\in\mc{N}$. On the other hand, $g(\tau)(F-p)(\tau)$ is holomorphic in $\HH$. Let $\mc{F}$ be a fundamental region for $\Gamma$ in $\HH$. Then we see that
\begin{eqnarray*}
\int_{\mc{F}} |g_a(\tau)|^2 (\tau-\bar{\tau})^k dudv &=& \int_{\mc{F}} \frac{\partial}{\partial \bar{\tau}} (g_a(\tau)\hat{G}_a(\tau)) dudv\\
 &=& \int_{\mc{F}} \frac{\partial}{\partial \bar{\tau}} (g_a(\tau)(F_a + \hat{G}_a - p_a)(\tau)) dudv.
\end{eqnarray*}
If we apply the Stoke theorem and take a sum, then we obtain that
\begin{equation} \label{applystoke}
\int_{\mc{F}} \sum_{a\in\mc{N}} |g_a(\tau)|^2 (\tau-\bar{\tau})^k dudv = \frac{-i}{2}  \int_{\partial\mc{F}}  \sum_{a\in\mc{N}} g_a(\tau)(F_a + \hat{G}_a - p_a)(\tau) d\tau.
\end{equation}
Since $\rho'$ is unitary, it follows that $\sum_{a\in\mc{N}} g_a(\tau)(F_a + \hat{G}_a - p_a)(\tau)$ is $\Gamma$-invariant. This implies that the right hand side of (\ref{applystoke}) is zero, since the sides of $\partial\mc{F}$ are paired by transformations in $\Gamma$. From this, we obtain that $g = 0$. Since $g=0$, by the definition of $\hat{G}(\tau)$ we see that $\hat{G} = 0$.  Then $(F-p)(\tau)$ is holomorphic on $\HH$ and at all cusps on $\Gamma$, from which we see that $F-p$ is zero since $-k<0$.  Thus, $f(\tau) = F^{k+1}(\tau) = 0$. By (\ref{Jacobithetaexpansion}) and (\ref{skewJacobithetaexpansion}), we see that $\Phi = \Psi = 0$. This completes the proof of injectivity.

Lastly, we prove the surjectivity of $\tilde{\eta}$. This proof is based on the argument in \cite{G}.
Suppose that $\{ g_{(\gamma,X)}|\ (\gamma,X)\in\Gamma^J\}$ is a parabolic cocycle in $P^e_{k,m}$. Then for every $(\gamma,X)\in\Gamma^J$ a function $g_{(\gamma,X)}(\tau,z)$ has the theta expansion
\[g_{(\gamma,X)}(\tau,z) = \sum_{a\in\mc{N}} g_{\gamma,a}(\tau)\theta_{2m,a,0}(\tau,z).\]
If we define $g_\gamma(\tau) = \sum_{a\in\mc{N}}g_{\gamma,a}(\tau)\mbf{e}_a$, then a collection $\{g_\gamma|\ \gamma\in\Gamma\}$ is a parabolic cocycle in $P_k$.
By Theorem \ref{converse} there exists a vector-valued function $H(\tau)$, holomorphic in $\HH$, such that
\[(H_0|_{-k,\chi',\rho'}\gamma)(\tau) - H_0(\tau) = g_\gamma(\tau)\]
for all $\gamma\in\Gamma$.
By \cite[Proposition 3.17]{CL2}, there exists $H_1(\tau)\in M^!_{-k,\chi',\rho'}(\Gamma)$ such that except possible finite term in the expansion $H_1(\tau)$ at the infinite cusp $H_1(\tau)$ has the same principal parts as $H_0(\tau)$. Thus, if we let $H(\tau) = H_0(\tau)-H_1(\tau)$, then $H(\tau)$ has poles only at the infinite cusp.
Let $W(\tau) = \frac{1}{c_{k+2}} D^{k+1}(H)(\tau)$, where $c_{k+2} = -\frac{(k-2)!}{(2\pi i)^{k-1}}$. Note that by the Bol's identity we have
\[c_{k+2}(W|_{k+2,\chi',\rho'}\gamma)(\tau) = D^{k+1}(H|_{-k,\chi',\rho'}\gamma)(\tau) = D^{k+1}(H+ g_\gamma)(\tau) = D^{k+1}(H)(\tau) = c_{k+2} W(\tau).\]
This implies that $W(\tau)\in M^!_{-k,\chi',\rho'}(\Gamma)$.
By (\ref{Fourierofvvpoincare}) we have a finite linear combination of vector-valued Poincar\'e series of weight $k+2$, multiplier system $\chi'$ and type $\rho'$, say $\sum_{l=1}^t b_lP(n_l,\alpha_l,\chi',\rho')(\tau)$, such that it has the same principal part with $W(\tau)$. Then we can write
\[W(\tau) = \sum_{l=1}^t b_lP(n_l,\alpha_l,\chi',\rho')(\tau) + f(\tau),\]
where $f(\tau)\in S_{k+2,\chi',\rho'}(\Gamma)$. Let $g(\tau) = \sum_{l=1}^t \overline{b_l}P_{n_l', \alpha_l, \overline{\chi'},\overline{\rho'}}(\tau)\in S_{k+2,\overline{\chi'},\overline{\rho'}}(\Gamma)$ (for the definition of $n_l'$ see section \ref{section2.3}). Then $g^*(\tau) = \sum_{l=1}^t b_lP(n_l,\alpha_l,\chi',\rho')(\tau)$. Therefore, we see that by \cite[Proposition 3.4]{CKL}
\[r(f,\gamma;\tau) + r(g^*,\gamma;\tau) = r(W,\gamma;\tau) = (H|_{-k,\chi',\rho'}\gamma)(\tau) - H(\tau) = g_{\gamma}(\tau).\]
If we define $\Phi(\tau,z) = \sum_{a\in\mc{N}} f_a(\tau)\theta_{2m,a,0}(\tau,z)$ and $\Psi(\tau,z) = \sum_{a\in\mc{N}} \overline{g_a(\tau)}\theta_{2m,a,0}(\tau,z)$, then we see that $\Phi(\tau,z)\in S_{k+2+\frac12, m,\chi}(\Gamma^J),\ \Psi(\tau,z) \in J^{sk,cusp}_{k+2+\frac12,m,\chi}(\Gamma^J)$ and
\[\tilde{\eta}(\Phi,\Psi) = \tilde{\beta}(\Phi) + \tilde{\alpha}(\Psi) = < g_{(\gamma,X)}|\ (\gamma,X)\in\Gamma^J>.\]
This proves the surjectivity of $\tilde{\eta}$, which completes the proof.

\bigskip

 

\begin{thebibliography}{99}



\bibitem{B} G. Bol, Invarianten linearer Differentialgleichungen, Abh. Math. Sem. Univ. Hamburg, 16: 3-4 (1949), 1-28.





\bibitem{BR} K. Bringmann and O. K. Richter, Zagier-type dualities and lifting maps for harmonic Maass-Jacobi forms, Adv. Math. 225 (2010), no. 4, 2298-2315.




\bibitem{CL3} D. Choi and S. Lim, Structures for pairs of mock modular forms with the Zagier duality, arXiv:1212.0370.

\bibitem{CL2} D. Choi and S. Lim, The Eichler cohomology theorem for Jacobi forms, arXiv:1211.2988.


\bibitem{CKL} D. Choi, B. Kim and S. Lim, Eichler integrals and harmonic weak Maass forms, arXiv:1210.3783.

\bibitem{CK} Y. Choie and H. Kim, An analogy of Bol's result on Jacobi forms and Siegel modular forms, J. Math. Anal. Appl. 257 (2001), no. 1, 79-88.




\bibitem{CL} Y. Choie and S. Lim, Eichler integrals, period relations and Jacobi forms, Math. Z. 271 (2012), no. 3-4, 639-661.



\bibitem{EO} T. Eguchi, H. Ooguri and Y. Tachikawa, Notes on the $K3$ surface and the Mathieu group $M24$, Exp. Math. 20, no. 1 (2011), 91-96.


\bibitem{E} M. Eichler, Eine Verallgemeinerung der Abelschen Integrale, Math. Z., 67 (1957), 267-298.






\bibitem{EZ} M. Eichler and D. Zagier, The theory of Jacobi forms, Birkhauser 1985.


\bibitem{G} J. Gim\'enez, Fourier coefficients of vector-valued modular forms of negative weight and Eichler cohomology, Thesis (Ph.D.)-Temple University, 2007, 87 pp.

\bibitem{GZ} L. G$\ddot{\mathrm{o}}$ttsche and D. Zagier, Jacobi forms and the structure of Donaldson invariants for $4$-mainfolds with $b_+ = 1$, Selecta Math. (N.S.) 4 (1) (1998) 69-115.


\bibitem{G2} R. Gunning, The Eichler cohomology groups and automorphic forms, Trans. Amer. Math. Soc. 100 (1961), 44-62.


\bibitem{HK} S. Husseini and M.  Knopp, Eichler cohomology and automorphic forms, Illinois J. Math. 15 (1971), 565-577.

\bibitem{Kon} M.  Knopp, Construction of automorphic forms on $H$-groups and supplementary Fourier series, Trans. Amer. Math. Soc. {\bf 103} (1962), 168--188.


\bibitem{K} M. Knopp, Some new results on the Eichler cohomology of automorphic forms, Bull. Amer. Math. Soc. 80 (1974), 607-632.


\bibitem{K2} M. Knopp, Modular Functions in Analytic Number Theory, 2nd Edition, Chelsea, New York, 1993.


\bibitem{KLR} M. Knopp, J. Lehner and W. Raji, Eichler cohomology for generalized modular forms, Int. J. Number Theory 5 (2009), no. 6, 1049-1059.




\bibitem{KM2} M. Knopp and H. Mawi, Eichler cohomology theorem for automorphic forms of small weights, Proc. Amer. Math. Soc. 138 (2010), no. 2, 395-404.




\bibitem{K2} I. Kra, On cohomology of Kleinian groups, Ann. of Math. (2) 89 (1969), 533-556.




\bibitem{Leh} J. Lehner, The Fourier coefficients of automorphic forms on horocyclic groups. II, Michigan Math. J. 6 (1959), 173-193.



\bibitem{L2} J. Lehner, The Eichler cohomology of a Kleinian group, Math. Ann., 192 (1971), 125-143.

\bibitem{L3} J. Lehner, Cohomology of vector-valued automorphic forms, Math. Ann., 204 (1973), 155-176.




\bibitem{Sko1} N. -P. Skoruppa, Developments in the theory of Jacobi forms, Acad. Sci. USSR, Inst. Appl. Math., Khabarovsk, 1990, 167-185.

\bibitem{Sko2} N. -P. Skoruppa, Explicit formulas for the Fourier coefficients of Jacobi and elliptic modular forms, Invent. Math. 102 (3) (1990) 501-520.



\bibitem{Zag} D. Zagier, Sur la conjecture de Saito-Kurokawa, in: Seminar on Number Theory, Paris 1979-80, in: Progr. Math., vol. 12, Birkh$\ddot{\mathrm{a}}$user, 1981, 371-394.

\bibitem{Z0} C. Ziegler, Jacobi forms of Higher Degree, Abh. Math. Sem. Univ. Hamberg 59, 191-224 (1989).


\bibitem{Z} S. Zwegers, Mock theta functions, PH.D Thesis, Universiteit Utrecht, 2002.
\end{thebibliography}
\end{document}